\documentclass[12pt]{article}
\usepackage{amsfonts}
\usepackage{amsmath}
\begin{document}

\newcommand{\li}{\underline}
\newcommand{\PP}{{\cal P}}
\newcommand{\RR}{{\cal R}}
\newcommand{\lam}{\lambda}
\newcommand{\ei}{\epsilon_i}
\newcommand{\be}{\bar{\epsilon}}
\newcommand{\eproof}{\hfill $\Box$ \vspace{0.2cm}}
\newcommand{\qed}{\hfill$\square$}

\newtheorem{theorem}{Theorem}
\newtheorem{proposition}{Proposition}
\newtheorem{rem}{Remark}
\newtheorem{corol}{Corollary}
\newtheorem{defi}{Definition}
\newtheorem{notat}{Notation}
\newtheorem{lemm}{Lemma}

\author{Nicholas J. Daras \thanks{Department of Mathematics, Faculty of Military Sciences, Hellenic Military Academy, 166 73, Vari Attikis, Greece,  E-mail: darasn@sse.gr} ~and Vassilios Nestoridis \thanks{Department of Mathematics, University of Athens, Panepistemiopolis, 157 84, Athens, Greece, E-mail: vnestor@math.uoa.gr}}

\title{UNIVERSAL TAYLOR SERIES ON CONVEX SUBSETS OF $\mathbb{C}^{n}$}
\date{}
\maketitle

\begin{abstract} We prove the existence of holomorphic functions $f$ defined on any open convex subset$\ \ {\rm \Omega}\subset {{\mathbb C}}^n$, whose partial sums of the Taylor developments approximate uniformly any complex polynomial on any convex compact set disjoint from $\overline{{\rm \Omega}}$ and on denumerably many convex compact sets in ${{\mathbb C}}^n\backslash {\rm \Omega}$ which may meet the boundary $\partial {\rm \Omega}$. If the universal approximation is only required on convex compact sets disjoint from$\ \ \overline{{\rm \Omega}}$, then $f$ may be chosen to be smooth on $\partial {\rm \Omega}$, that is $f\in A^{\infty }\left({\rm \Omega}\right)$. Those are generic universalities.\end{abstract}

\noindent \textbf{Subject Classification MSC2010: }primary 30K05, 32A05, 32A30, secondary 40A05, 41A58, 41A99\textbf{}


\section{ Introduction }

Analysis is the study of limiting processes. Not every limiting process, of course, converges, but examples have been found where processes diverge in a maximal way. Such an extreme behavior is often linked with the phenomenon of \textit{universality}. The present paper lies with multidimensional universality. 

The first universal series was observed by Fekete before 1914. He showed that there exists a formal real power series on [-1\textit{; }1] that not only diverges at every point different from $0$ but does so in the worst possible way. In 1951, Seleznev proved the existence of a (divergent) power series with universal approximation properties in ${\mathbb C} \backslash \left\{0\right\}$. In 1970 and 1971, Luh and Chui and Parnes showed the existence of a universal Taylor series with non zero radius of convergence and universal properties outside the closure of the disk of convergence of the Taylor series (\cite{Chui}, \cite{Luh}). In 1987, Grosse-Erdmann proved that all these universalities (and many others) were generic and, thus, the use of Baire's theorem could yield considerable simplifications in the proofs of the already known universalities (\cite{Gros1}). Improved forms of universality results have been obtained by Nestoridis (\cite{Nest1}, \cite{Nest2} and \cite{Nest3}), Melas and Nestroridis (\cite{Mel1} and \cite{Mel2}), M$\ddot{u}$ller, V. Vlachou and A. Yavrian (\cite{Mul}), and Nestroridis and Papadimitropoulos (\cite{Nest4}). Now, there is a systematic use of Baire's theorem in order to establish new generic universalities.

The theory of universal Taylor series in several complex variables is difficult to be developed because of lack of general approximation theorems as Mergelyan's or Runge's theorem. Recently, Clou$\hat{a}$tre (\cite{Clo}) obtained Seleznev-type universal Taylor series in several complex variables. Now, the universal approximation is valid on any polynomially convex compact subset $K\subset \subset{{\mathbb C}}^n$ that is disjoint from the origin. The Taylor series in \cite{Clo} do not converge on any non-empty open subset of$\ \ {{\mathbb C}}^n$. 

In the present article, we prove the existence of holomorphic functions on any open convex subset$\ \ {\rm \Omega}\subset{{\mathbb C}}^n$, such that the partial sums of their Taylor development approximate any analytic polynomial $P\left(z\right)$ on any convex compact set disjoint from $\overline{{\rm \Omega}}$ and on countably many convex sets contained in ${{\mathbb C}}^n\backslash{\rm \Omega}$ and possibly meeting the boundary$\ \ \partial {\rm \Omega}$ of$\ \ {\rm \Omega}$. If we restrict our attention only on convex compact sets disjoint from$\ \ \overline{{\rm \Omega}}$, then the universal function $f$ may be chosen to be smooth on$\ \partial {\rm \Omega}$, that is $ f \in A^{\infty} \left({\rm \Omega}\right)$.

The key properties which make possible our proofs are Kallin's lemma and the Oka-Weil Theorem.\\
\\
\noindent \textbf{Lemma 1.1} (Kallin's Lemma; see, for instance, \cite{Kal} and \cite{Pae}) \textit{Suppose }$K_1$\textit{ and }$K_2$\textit{ are polynomially convex subsets of }${{\mathbb C}}^n$\textit{, suppose there is a polynomial }$p$\textit{ mapping }$K_1$\textit{ and }$K_2$\textit{ onto two polynomially convex subsets }$F_1$\textit{ and }$F_2$\textit{ of the complex plane such that }$0$\textit{ is a boundary point of both }$F_1$\textit{ and }$F_2$\textit{ and with }$F_1\bigcap F_2=\left\{0\right\}$\textit{. If }$p^{-1}\left(0\right)\bigcap \left(K_1\cup K_2\right)$\textit{ is polynomially convex, then }$K_1\cup K_2$\textit{is polynomially convex. }\\

According to Kallin's lemma, \textit{if} $K_1$ \textit{and} $K_2$ \textit{are two compact sets in} ${{\mathbb C}}^n$ \textit{and} $p$ \textit{a polynomial such that the polynomial convex hulls} $\widehat{p\left(K_1\right)}$ and $\widehat{p\left(K_2\right)}$ \textit{are disjoint, then the polynomial convex hull }$\widehat{K_1 \cup K_2}$ \textit{of the union }$K_1\cup K_2$\textit{ equals the union} $\widehat{K_1}\cup\widehat{K_2}$ \textit{of the polynomial convex hulls}. Using Kallin's lemma and the fact that \textit{every compact convex set in }${{\mathbb C}}^n$\textit{ is polynomially convex}, one deduces that \\
\\
\noindent \textbf{Lemma 1.2}\textit{ The union }$\ K_1\cup K_2\ $\textit{ of two disjoint compact convex sets }$K_1$ \textit{and} $K_2$ \textit{is polynomially convex}. \\

 Then, application of the Oka-Weil theorem (see, for instance, \cite{Ran}) guarantees that \\
\\
\noindent \textbf{Lemma 1.3}\textit{ Any holomorphic function }$h$\textit{ on an open set} $V\supset K_1\cup K_2$ \textit{can be uniformly approximated uniformly on} $K_1\cup K_2$ \textit{by holomorphic polynomials, provided that }$K_1$ \textit{and} $K_2$ are two disjoint compact convex subsets of ${{\mathbb C}}^n$. \\

Further, again by Kallin's lemma \\
\\
\noindent \textbf{Lemma 1.4}\textit{ If }${{\rm \Omega}}_1$ \textit{and} ${{\rm \Omega}}_2$ \textit{are two disjoint open convex subsets of}$\ {{\mathbb C}}^n$\textit{, then}$\ \ {{\rm \Omega}}_1 \cup {{\rm \Omega}}_2$ \textit{is a Runge domain}.\\

 It follows that \textit{every holomorphic function} $f$ \textit{on} ${{\rm \Omega}}_1 \cup {{\rm \Omega}}_2$ \textit{may be approximated uniformly on compact sets by a sequence of holomorphic polynomials}. Hence, by Weierstrass theorem, we obtain the following.\\
\\
\noindent \textbf{Lemma 1.5}\textit{ There exists a polynomial }$Q$ \textit{approximating} $f$ \textit{uniformly on any prescribed compact set} $L\subset{{\rm \Omega}}_1 \cup {{\rm \Omega}}_2$ \textit{and such that a finite set of partial derivatives} $D^{\left(\ell \right)}Q$ \textit{approximate uniformly on} $L$ \textit{the corresponding partial derivatives} $D^{\left(\ell \right)}f$ of$\ \ f$, \textit{provided that }${{\rm \Omega}}_1$ \textit{and} ${{\rm \Omega}}_2$ \textit{are two disjoint open convex subsets of} ${{\mathbb C}}^n$ \textit{and }$f\ $\textit{is holomorphic in} ${{\rm \Omega}}_1\cup {{\rm \Omega}}_2$.\\

In this paper, the above facts from the theory of functions of several complex variables are combined with the abstract theory of universal series (\cite{Bay}, \cite{Nest4}). To do this, we fix any enumeration$\ \ {\left(N_j\right)}_{j=0,1,2,\dots }$ of ${{\mathbb N}}^n$ and we order all monomials in the Taylor development of any holomorphic function using this enumeration. The universal approximation is carried out by partial sums of the form
\[S_\lambda\left(f,\zeta\right)\left(z\right)=\sum^\lambda_{j=0}{a_{N_j}\left(f,\zeta\right){\left(z-\zeta\right)}^{N_j}},\] 
where

\noindent $\ \ z=\left(z_1,z_2,\dots ,z_n\right)$, $\zeta=\left(\zeta_1,\zeta_2,\dots ,\zeta_n\right)$, $N_j=\left(N_{j,1},N_{j,2},\dots ,N_{j,n}\right)$,\\
 ${\left(z-\zeta\right)}^{N_j}={\left(z_1-\zeta_1\right)}^{N_{j,1}}\dots {\left(z_n-\zeta_n\right)}^{N_{j,n}}$ and 
\[a_{N_j}\left(f,\zeta\right):=\frac{1}{\left(N_{j,1}\right)!\dots \left(N_{j,n}\right)!}\ \frac{{\partial }^{N_{j,1}+\dots +N_{j,n}}}{\partial z^{N_{j,1}}_1\dots \partial z^{N_{j,n}}_n}f\left(\zeta\right).\] 

The abstract theory of universal series allows that, if we start by any infinite set$\ \ {\mathbb M}\subseteq {\mathbb N}$, then we can always require that all the indices$\ \ \lambda_N$, such that the corresponding partial sum $S_{\lambda_N}\left(f,\zeta\right)\left(z\right)$ do the universal approximation, may be chosen in$\ {\mathbb M}$; that is $\lambda_N\in{\mathbb M}$ for all$\ N=1,2,\dots $. Therefore, choosing the enumeration$\ \ {\left(N_j\right)}_{j=0,1,2,\dots }$ of ${{\mathbb N}}^n$ properly and the infinite set ${\mathbb M}\subseteq{\mathbb N}$ in a good way, we see that the partial sums realizing the approximation can take the form 
\[S_{\lambda_N}\left(f,\zeta\right)\left(z\right)=\sum_{\nu_1+\nu_2+\dots +\nu_n\leq \tau_N}{a_{\nu_1,\nu_2,\dots ,\nu_n}\left(f,\zeta\right){\left(z_1-\zeta_1\right)}^{\nu_1}\dots {\left(z_n-\zeta_n\right)}^{\nu_n}}\] 
or the form 
\[S_{\lambda_N}\left(f,\zeta\right)\left(z\right)=\sum_{\nu^2_1+\nu^2_2+\dots +\nu^2_n\leq\tau^2_N}{a_{\nu_1,\nu_2,\dots ,\nu_n}\left(f,\zeta\right){\left(z_1-\zeta_1\right)}^{\nu_1}\dots {\left(z_n-\zeta_n\right)}^{\nu_n}}. \] 
Certainly the results given below are much more general. 

The abstract theory of universal series, which is used in the present article, is based on Baire's category theorem. For the role of Baire's theorem in Mathematical Analysis we refer to \cite{Gros2} and \cite{Kah}.

Of course, a natural question which may now be asked is the following. Can our universality results generalize to more general domains and to more general compact sets in$\ {{\mathbb C}}^n$? The general answer depends strongly on the polynomial convexity of a finite union of disjoint polynomially convex compact subsets of${\rm \ }{{\mathbb C}}^n$. In this direction, we recall some basic and elementary facts concerning unions of polynomially convex compact sets. 

\begin{itemize}
\item[(i)] \textit{The union of a polynomially convex compact set and a finite set of points is polynomially convex }(see, for instance Remark 4 in \cite{Pet}).
\item[(ii)] \textit{If }$K_1\cup \ K_2$\textit{ is polynomially convex and compact, }$K_1\bigcap K_2=\emptyset $\textit{ and }$K'_1\subset K_1$\textit{ is polynomially convex and compact, then }$K'_1\cup K_2$\textit{ is polynomially convex} (see, for instance Remark 4 in \cite{Pet}). 
\item[(iii)] \textit{The union of a polynomially convex compact set and finitely many disjoint smooth compact curves in its complement is polynomially convex} (\cite{Sto}). Notice that the smoothness of the compact curves is an essential condition. For instance, there exist (non-smooth) arcs with non-trivial polynomial convex hulls: J. Wermer gave the first example of an arc in ${{\mathbb C}}^{{\rm 3}}$ that is not polynomially convex (\cite{Wer}). Furthermore, the notion of polynomial convexity is not invariant under biholomorphic mappings. This phenomenon was first observed by J.Wermer (\cite{Ex}).
\item[(iv)] \textit{Let }$K$\textit{ be a compact set disjoint from the closed ball }$\overline{B}$. \textit{Assume that }$K$\textit{ is a finite union of disjoint polynomially convex sets}. \textit{Then for every }$\epsilon\ >\ 0$ \textit{there is a compact set }$K_\epsilon\subset K$, \textit{such that }$\overline{B}\ \cup K_\epsilon\ $\textit{ is polynomially convex and }$\left|K\ \backslash \ K_\epsilon\right|<\epsilon$\textit{ }(see, for instance, Lemma 3.6 in \cite{For}).
\end{itemize}

Further, for three polynomially convex compact sets in ${{\mathbb C}}^n$, we also have the next promising result. \textit{If }$B_1$\textit{,}$\ B_2$\textit{ and }$B_3$\textit{ are three pairwise disjoint closed balls in }${{\mathbb C}}^n$\textit{ then their union }$\ B_1\cup \ B_2\cup B_3$\textit{ is polynomially convex }(\cite{Kal2}). However, \textit{there exist three congruent, pairwise disjoint, closed polydisks }$P_1$\textit{,}$\ P_2$\textit{ and }$P_3$\textit{ in }${{\mathbb C}}^3$\textit{ such that }$P_1\cup \ P_2\cup P_3$\textit{ is not polynomially convex }(\cite{Kal2}). Notice that Kallin's proof actually used polydisks parallel to the coordinate axes. This is, however, not possible in ${{\mathbb C}}^2$ (see \cite{Ros}). Another counterexample was proved in 1984 by A. M. Kytmanov and G. Khud$\check{a}$iberganov: \textit{there exist three congruent, pairwise disjoint, closed complex ellipsoids}$\ \ E_1$\textit{,}$\ E_2$\textit{ and }$E_3$\textit{ in }${{\mathbb C}}^3$\textit{ such that }$E_1\cup \ E_2\cup E_3$\textit{ is not polynomially convex }(\cite{Khu}, \cite{Kyt}). 

The first example of three pairwise disjoint compact convex sets in ${{\mathbb C}}^2$ whose union is not polynomially convex was published by J.-P. Rosay in 1989: \textit{there exist three congruent, pairwise disjoint, convex closed limited tubes }$T_1$\textit{,}$\ T_2$\textit{ and }$T_3$\textit{ in }${{\mathbb C}}^2$\textit{ such that }$T_1\cup \ T_2\cup T_3$\textit{ is not polynomially convex}. Here the limited tube in ${{\mathbb C}}^2$ with base domain $B\ \subset {{\mathbb R}}^2$ and height $M$ is the domain $\{(z_1,\ z_2)\in{{\mathbb C}}^2\ :\ (Re\ z_1,\ Re\ z_2)\in\ B,\ |Im\ z_1|\ <\ M,\ |Im\ z_2|\ <\ M\}$ (\cite{Ros}). 

Finally, in 1996, U. Backlund and A. F$\ddot{a}$llstr$\ddot{o}$m proved the last well known counterexample: \textit{There exist a positive integer }$k$\textit{ and three pairwise disjoint, closed sets }$S_1$\textit{,}$\ S_2$\textit{ and }$S_3$\textit{ in }${{\mathbb C}}^3$\textit{ all congruent to }$\left\{(z_1, z_2, z_3)\ \in {{\mathbb C}}^3\ :\ {|z_1|}^2 + {|z_2|}^2 + {|z_3|}^{2k}\ \leq \ 1\right\}$ \textit{such that }$S_1\cup \ S_2\cup S_3$\textit{ is not polynomially convex} (\cite{Back}). Domains of the form 
 $\{(z_1, z_2, z_3) \in {{\mathbb C}}^3 :\ {|z_1|}^2 + {|z_2|}^2 + {|z_3|}^{2k} \leq  1\}$ have been studied by E. Bedford and S. Pinchuk (\cite{Bed}). 

Summarizing, the union of three convex disjoint compact subsets of ${{\mathbb C}}^n$ is not polynomially convex in general. This is an obstacle in order to obtain universal approximation on more general compact sets which are not convex. However, in the last section of the present paper we state some weak universality results where the universal approximation is valid on finite unions of disjoint smooth curves. This is obtained for special domains ${\rm \Omega}\subset {{\mathbb C}}^n$ as domains of definition of the universal functions not necessarily convex.

\section{Universalities in ${\mathcal O}\left({\mathbf \Omega}\right)$}

Let ${\rm \Omega}\in{{\mathbb C}}^n$ be an open convex set and let ${\mathcal O}\left({\rm \Omega}\right)$ the space of all holomorphic functions in$\ {\rm \Omega}$. Let also $\zeta\in{\rm \Omega}$, $\zeta=\left(\zeta_1,\zeta_2,\dots ,\zeta_n\right)$, be a fixed point in ${\rm \ \Omega}$. Every holomorphic function $f\in{\mathcal O}\left({\rm \Omega}\right)$ has a Taylor development:
\[f\left(z\right)=\sum_{\nu\in{{\mathbb N}}^n}{a_\nu\left(f,\zeta\right){\left(z-\zeta\right)}^\nu}\] 
locally around $\zeta$, where $\nu=\left(\nu_1,\nu_2,\dots ,\nu_n\right)\in{{\mathbb N}}^n$ and ${\left(z-\zeta\right)}^\nu={\left(z_1-\zeta_1\right)}^{\nu_1}\dots {\left(z_n-\zeta_n\right)}^{\nu_n}$. Let $N_j$, $j=0,1,2,\dots $ be an enumeration of ${{\mathbb N}}^n$. Then $f\in{\mathcal O}\left({\rm \Omega}\right)$ can be identified with the sequence
\[\left(a_{N_j}\left(f,\zeta\right):j=0,1,2,\dots \right).\] 

Let 
\[{\mathcal A}:=\left\{{\left(a_{N_j}\left(f,\zeta\right)\right)}^\infty_{j=0}:f\in{\mathcal O}\left({\rm \Omega}\right)\ \right\}. \] 
The space ${\mathcal O}\left({\rm \Omega}\right)$ endowed with the topology of uniform convergence on compact subsets of ${\rm \Omega}$ is a Fr$\acute{e}$chet space. A standard metric on ${\mathcal O}\left({\rm \Omega}\right)$ yields naturally a metric$\ \ d\ \ $on$\ \ {\mathcal A}$, such that the map
\[{\mathcal O}\left({\rm \Omega}\right)\ni f\mapsto{\left(a_{N_j}\left(f,\zeta\right)\right)}^\infty_{j=0}\in{\mathcal A}\] 
becomes a surjective isometry. Thus, $\left({\mathcal A},d\ \right)$ is a complete space. It is easy to see that \\

\noindent \textbf{Proposition 2.1 }
\begin{itemize}
\item[1)]The projections
\[{\mathcal A}\ni{a=\left(a_j\right)}^\infty_{j=0}\mapsto a_k\in{\mathbb C}\] 
are continuous for every $k=0,1,2,...$ 
\item[2)] The set $C_{\rm oo}:=\left\{a=\left(a_0,a_1,\dots ,a_M,0,0,0,\dots \ \ \right):M\ge 0\right\}$ is contained in ${\mathcal A}$ 
\item[3)]  The set $C_{\rm oo}$ is dense in ${\mathcal A}$ (with respect to the metric $d$). 
\end{itemize}
The last statement is true, since, by the Oka-Weil theorem, the set of polynomials is dense in ${\mathcal O}\left({\rm \Omega}\right)$, for any convex open subset ${\rm \Omega}\in{\mathbb C}^n$ (see, for instance, \cite{Ran}).

Now let $K\in{\mathbb C}^n$ be a compact convex set, such that$\ \ K \cap {\rm \Omega}=\emptyset $. Let also $\mathfrak{X}$ be the space of restrictions of all polynomials on $K$ endowed with the metric $\varrho$ induced by the supremum norm on$\ \ K$. Let finally$\ \ \ x_j\in\mathfrak{X}$, $j=0,1,2,...$ denote the sequence$\ x_j={\left(z-\zeta\right)}^{N_j}$, where $N_j\ $is the enumeration of$\ \ {\mathbb N}^n$ considered previously ($j=0,1,2,...$). Let also ${\mathbb M}$\textit{ }be an infinite subset of\textit{ }${\mathbb N}$.\\
\\ 
\textbf{Definition 2.2} \textit{We say that a numerical sequence}
$$a=\left(a_j\right)\in{\mathcal A}$$
\textit{  belongs to the class} $${\mathcal U}^{\left({\mathbb M},\left(N_j\right)\right)}_{\mathcal A}$$
\textit{ of universal multiple series (with respect to the enumeration }$N_j$\textit{ of }${\mathbb N}^n$\textit{ and the infinite set }${\mathbb M}\subset {\mathbb N}$\textit{), if for every }$P\in \mathfrak{X}$\textit{ there exists a sequence } 
$${\left(\lambda_m\in {\mathbb M}\right)}_{m=0,1,2,\dots \ }$$
\textit{satisfying the following. }
\begin{itemize} 
\item[i.] $\sum^{\lambda_{m}}_{j=0}{a_j x_j} \underset{m\rightarrow \infty}{\longrightarrow} P$ \textit{ in }$\mathfrak{X}$\textit{, with respect to }$\varrho$.
\item[ii.] \textit{$\sum^{\lambda_m}_{j=0}{a_{n_j}e_j}\underset{m\rightarrow \infty}{\longrightarrow}a$}\textit{, with respect to }$d$\textit{.} 
\end{itemize}
\textit{where} $e_j=\left(0,\dots ,0,1,0,\dots \right)$ \textit{and}$\ $1 \textit{appears in the} $j$\textit{${}^{th}$} \textit{place}. \\

According to \cite{Bay}, \cite{Nest4}, a necessary and sufficient condition to hold
\[{{\mathcal U}}^{\left({\mathbb M},\left(N_j\right)\right)}_{{\mathcal A}}\neq \emptyset \] 
is that, for every $P\in \mathfrak{X}$ and every $\epsilon>0$, there exists   $ \ b \ =\left(b_0,b_1,\dots ,b_M,0,0,0,\dots \ \ \right)\in C_{{\rm oo}}$ such that 

$$d\left(\sum^M_{j=0}{b_je_j},\ 0\right)<\epsilon$$ 
and  
$$\varrho\left(\sum^M_{j=0}{b_jx_j},\ P\right)<\epsilon.$$

\noindent Translating to our context, it suffices to show that, for each compact convex set $L\subset \subset{\rm \Omega}$ and any polynomial $P\left(z\right)\in \mathfrak{X}$, there exists a complex polynomial $Q\left(z\right)$ such that 

\begin{center}
 ${sup}_{z\in L}\left|Q\left(z\right)\right|<\epsilon$ and ${sup}_{z\in K}\left|Q\left(z\right)-P\left(z\right)\right|<\epsilon$.
\end{center}

 To do so, let us consider the compact set $L\cup K$. By Lemma 1.2, this is a polynomially convex compact set, as the union of two disjoint compact convex subsets of ${{\mathbb C}}^n$. The function 
$$w:L\cup K \to {\mathbb C}:z\mapsto w\left(z\right)=\left\{ \begin{array}{c}
0,\ on\ L \\ 
P\left(z\right),\ on\ K \end{array}
\right.$$
is holomorphic in a neighbourhood of $L\cup K$. From the Oka-Weil theorem, it follows that there exists a polynomial $Q\left(z\right)$ such that
\[{sup}_{z\in L\cup K}\left|w\left(z\right)-Q\left(z\right)\right|<\epsilon \ \ (\cite{Ran}).\] 
Thus, ${sup}_{z\in L}\left|Q\left(z\right)\right|<\epsilon$ and ${sup}_{z\in K}\left|Q\left(z\right)-P\left(z\right)\right|<\epsilon$. We infer
\[{{\mathcal U}}^{\left({\mathbb M},\left(N_j\right)\right)}_{{\mathcal A}}\neq \emptyset {\rm \ }\ .\] 
Hence, according to \cite{Bay} or \cite{Nest4}, we have proved the following result.\\
\\
\textbf{Proposition 2.3} \textit{Let }${\mathbb M}$\textit{ be an infinite subset of }${\mathbb N}$\textit{. Let }${\rm \Omega}\subset {{\mathbb C}}^n$\textit{ be an open convex set,
}$\zeta\in{\rm \Omega}$\textit{ be fixed and }$K\subset {{\mathbb C}}^n$\textit{ be a compact convex set disjoint from }${\rm \Omega}$\textit{. Let also }${\left(N_j\right)}_{j=0,1,2,\dots
}$ \textit{be any enumeration of} ${{\mathbb N}}^n$\textit{. Then there exists a holomorphic function }$f\in{\mathcal O}\left({\rm \Omega}\right)$\textit{ whose partial sums }
$$ S_\lambda\left(f,\zeta\right)\left(z\right)=\sum^\lambda_{j=0}{a_{N_j}\left(f,\zeta\right){\left(z-\zeta\right)}^{N_j}}, \ \lambda=0,1,2,\dots $$
\textit{of the Taylor development of }$f$\textit{ around }$\zeta$\textit{ satisfy the following.}
\newline 

For every holomorphic polynomial $P\left(z\right)$\textit{, }$z\in{\mathbb C}^n$\textit{, there exists a sequence }${\left(\lambda_N\in{\mathbb M}\right)}_{N=1,2,\dots }$\textit{ such that}
\begin{center}
$S_{\lambda_N}\left(f,\zeta\right)\left(z\right)\rightarrow P\left(z\right)$
\textit{uniformly on} $K$ \textit{, as } $N \rightarrow \infty$
\end{center}
\textit{and}
\begin{center}
$S_{\lambda_N}\left(f,\zeta\right)\left(z\right)\rightarrow f\left(z\right)$\textit{ uniformly on each compact subset of }${\rm \Omega}$\textit{, as }$N \rightarrow \infty. $
\end{center}
The set of all such functions $f$\textit{ is a dense }${\mathcal G}_\delta$\textit{ subset of}$\ {\mathcal O}\left({\rm \Omega}\right)$ \textit{and contains a dense vector space except}
$\ 0$, \textit{where} ${\mathcal O}\left({\rm \Omega}\right)$ \textit{is endowed with the topology of uniform convergence on
 compacta}.\\
\\
\textbf{Remark 2.4 }Since, for any compact convex set $K\subset \subset {{\mathbb C}}^n$, the set of polynomials is dense in ${\mathcal O}\left(K\right)$, it follows that in Proposition 2.3 \textit{the partial sums} $S_\lambda\left(f,\zeta\right)\left(z\right)$ \textit{approximate uniformly on} $K$ \textit{any function holomorphic in any} (\textit{varying}) \textit{neighbourhood of} $K$. This applies to all the results of the present paper.\\

Combining Proposition 2.3 with Theorem 3 in \cite{Bay}, we can obtain a similar result simultaneously for all sets $K$ in any denumerable family of compact convex subsets of ${{\mathbb C}}^n$ disjoint from ${\rm \Omega}$. In particular, if$\ \ K$ is a compact convex set disjoint from$\ \overline{{\rm \Omega}}$, then it can be separated from $\overline{{\rm \Omega}}$ using a hyperplane defined by a real linear functional with coefficients in ${\mathbb Q}$. Thus, there exists a denumerable family of compact convex sets disjoint from $\overline{{\rm \Omega}}$ containing all other such $K$'s. Hence, we get the following result.\\
\\
\textbf{Theorem 2.5} \textit{Let }${\mathbb M}$\textit{ be an infinite subset of }${\mathbb N}$\textit{. Let }${\rm \Omega}\subset{{\mathbb C}}^n$\textit{ be an open convex set and }${\left(K_m\right)}_{m=1,2,\dots }$ be \textit{a family of compact convex subsets of} ${{\mathbb C}}^n$ \textit{disjoint with} ${\rm \Omega}$. \textit{Let }$\zeta\in{\rm \Omega}$\textit{ be fixed and let }${\left(N_j\right)}_{j=0,1,2,\dots }$ \textit{be any enumeration of} ${{\mathbb N}}^n$\textit{. Then there exists a holomorphic function }$f\in {\mathcal O}\left({\rm \Omega}\right)$\textit{ such that the partial sums} 
$$S_\lambda\left(f,\zeta\right)\left(z\right)= \sum^\lambda_{j=0}{a_{N_j}\left(f,\zeta\right){\left(z-\zeta\right)}^{N_j}}, \ \lambda=0,1,2,\dots $$
\textit{of the Taylor development of }$f$\textit{ around }$\zeta$\textit{ satisfy the following.}\\
\textit{For every compact convex set} $K\subset \subset {{\mathbb C}}^n$\textit{ which is disjoint from }$\overline{{\rm \Omega}}$\textit{ or it is equal to }$K_m$\textit{ for some }$m=1,2,\dots $\textit{ and for every analytic polynomial }$P\left(z\right)$\textit{, }$z\in{{\mathbb C}}^n$\textit{, there exists a sequence }$\left(\lambda_N\in{\mathbb M}\ :N=1,2,\dots \right)$\textit{ such that} 
\begin{center}
$S_{\lambda_N}\left(f,\zeta\right)\left(z\right)\rightarrow P\left(z\right)$\textit{ uniformly on }$K$\textit{, as }$N\rightarrow \infty$
\end{center}
\textit{and}
\begin{center}
$S_{\lambda_N}\left(f,\zeta\right)\left(z\right)\rightarrow f\left(z\right)$ \textit{ uniformly on each compact subset of }${\rm \Omega}$\textit{, as }$N\rightarrow \infty.$
\end{center}
\textit{ The set of all such functions }$f$\textit{ is a dense }${{\mathcal G}}_\delta$\textit{ subset of}$\ \ {\mathcal O}\left({\rm \Omega}\right)$ \textit{and contains a dense vector space except}$\ 0$\textit{, where} ${\mathcal O}\left({\rm \Omega}\right)$ \textit{is endowed with the topology of uniform convergence on compacta}. \\
\\
\textbf{Remark 2.6 }If we allow the center of expansion $\zeta$ to vary in ${\rm \Omega}$, then, combining Theorem 2.5 with Theorem 3 in \cite{Bay}, we obtain the following stronger version of Theorem 2.5. For every compact subset $L\subset \subset {\rm \Omega}$ we have
$$ {sup}_{\zeta\in L,z\in K}\left|S_{\lambda_N}\left(f,\zeta\right)\left(z\right)-P\left(z\right)\right|{{\rightarrow 0}} $$
and 
$${sup}_{\zeta\in L,z\in L}\left|S_{\lambda_N}\left(f,\zeta\right)\left(z\right)-f\left(z\right)\right|{{\rightarrow\ 0}}$$ 
as $N\rightarrow \infty$. This result is also generic in the space ${\mathcal O}\left({\rm \Omega}\right)$, provided that the open set ${\rm \Omega}\subset {{\mathbb C}}^n$ is convex.\\
\\
\noindent \textbf{Remark 2.7 }For particular open convex sets ${\rm \Omega}\subset {{\mathbb C}}^n$, the result of Theorem 2.5 remain valid simultaneously for all compact convex sets $K$ disjoint from ${\rm \Omega}$. This is the case for instance if 
$${\rm \Omega}=\left\{z=\left(z_1,\dots ,z_n\right)\in {{\mathbb C}}^n:Rez_1>0\right\}. $$ 

\section{Universalities in ${{\mathbf A}}^{\infty }\left({\mathbf \Omega}\right)$}

Let ${\rm \Omega}\subset{{\mathbb C}}^n$ be an open convex set. A holomorphic function $f\in {\mathcal O}\left({\rm \Omega}\right)$ belongs to $A^{\infty }\left({\rm \Omega}\right)$, if every order's partial derivative of $f$ extends continuously on $\overline{{\rm \Omega}}$. 

We endow $A^{\infty }\left({\rm \Omega}\right)$ with the topology of uniform convergence of any order's partial derivatives on each compact subset of $\overline{{\rm \Omega}}$. It is well known that $A^{\infty }\left({\rm \Omega}\right)$\textit{ }is a $Fr\acute{e}chet$ space.\\
\\
\textbf{Proposition 3.1} \textit{The set of polynomials is dense in} $A^{\infty }\left({\rm \Omega}\right)$,\textit{ for every convex open set }${\rm \Omega}\subset {{\mathbb C}}^n$\textit{. }\\ 
\\
\textit{Proof} Without loss of generality, we may assume that $0\in {\rm \Omega}$. If $f\in A^{\infty }\left({\rm \Omega}\right)$ and $0<r<1$, we denote by $f_r$ the function $f_r\left(z\right)=f\left(rz\right)$ which is holomorphic in an open neighbourhood of $\overline{{\rm \Omega}}$. We can easily cheek that $f_r\rightarrow f$ in $A^\infty\left({\rm \Omega}\right)$ as $r\rightarrow 1$. Thus, $f$ can be approximated in the topology of $A^\infty\left({\rm \Omega}\right)$ by some $f_r$, $0<r<1$. Since $f_r$ is defined in an open neighbourhood of $\overline{{\rm \Omega}}$ and ${\rm \Omega}$ is convex, it follows easily that $f_r$ may be approximated in the topology of $A^\infty \left({\rm \Omega}\right)$ by polynomials \footnote{   A compact polynomially  convex  set $K \subset {{\mathbb C}}^n$  admits  a  basis  of  pseudoconvex   neighbourhoods     ${\left({{\rm \Omega}}_\nu\right)}_{\nu\in {\mathbb N}}$  that are Runge in   ${{\mathbb C}}^n$, in the sense that the algebra of polynomials   ${\mathcal P}$ is dense in every   ${\mathcal O}\left({{\rm \Omega}}_\nu\right).$ }(see, for instance, Theorem 4.12.1, page 143, in \cite{Forst}). Thus, the set of polynomials is dense in $A^\infty\left({\rm \Omega}\right)$, for any ${\rm \Omega}\subset {{\mathbb C}}^n$ being an open convex set.

According to \cite{Bay}, in order to establish existence of universal Taylor series in $A^{\infty }\left({\rm \Omega}\right)$, it suffices to prove the following.\\
 
\textbf{Lemma 3.2} \textit{Let }${\rm \Omega}\subset {{\mathbb C}}^n \ $ \textit{ be an open convex set. Let }$K\subset \subset {{\mathbb C}}^n$\textit{be a compact convex set disjoint from }$\overline{{\rm \Omega}}$\textit{ and let }$P\left(z\right)$\textit{ be a holomorphic polynomial on }${{\mathbb C}}^n$\textit{. Let also }$0<\delta,\ \epsilon<\infty $\textit{. Let }$\mathfrak{F}$\textit{ be a finite set of symbols of partial derivation. Then there exists a holomorphic polynomial }$Q\left(z\right)$\textit{ such that }
\begin{center}
${sup}_{z\in K}\left|Q\left(z\right)-P\left(z\right)\right|<\epsilon$\textit{ and }${sup}_{z\in \overline{{\rm \Omega}}\bigcap \overline{D\left(0;\delta\right)}}\left|D^{\left(\ell \right)}Q\left(z\right)\right|<\epsilon$\textit{for all }$\ell \in \mathfrak{F}.$
\end{center}
\textit{Here }\\
\\
\vspace{5 mm} $\bullet \overline{D\left(0;\delta\right)}=\left\{z\in {{\mathbb C}}^n:\left\|z\right\|\le \delta\right\}$\textit{ and}\\
\vspace{5 mm} $\bullet \left(D^{\left(\ell \right)}Q\right)$ \textit{ and } $\left(D^{\left(\ell \right)}P\right)$\textit{ denote partial derivatives of }$Q$\textit{ and }$P$\textit{ respectively, that come from the same symbol}
$$D^{\left(\ell \right)}=\left({{\partial }^{{\ell }_1+\dots +{\ell }_n}}/{\partial z^{{\ell }_1}_1\dots }\partial z^{{\ell }_n}_n\right), \ \ell =\left({\ell }_1,\dots ,{\ell }_n\right).$$

\textit{Proof} We remark that, since $K$ and $\overline{{\rm \Omega}}\bigcap \overline{D\left(0;\delta\right)}$ are disjoint compact convex sets in ${{\mathbb C}}^n$, we can find disjoint open convex sets ${{\rm \Omega}}_1$, ${{\rm \Omega}}_2\subset {{\mathbb C}}^n$ such that $\overline{{\rm \Omega}}\cup \overline{D\left(0;\delta\right)}\subset {{\rm \Omega}}_1$ and $K\subset {{\rm \Omega}}_2$. It is easily seen that ${{\rm \Omega}}_1\cup {{\rm \Omega}}_2$ is a Runge domain. Let us consider the holomorphic function
\[w:{{\rm \Omega}}_1\cup {{\rm \Omega}}_2\to {\mathbb C}{\rm :}w\left(z\right):=\left\{ \begin{array}{c}
0,\ on\ {{\rm \Omega}}_1 \\ 
P\left(z\right),\ on{{\rm \ \Omega}}_2.\  \end{array}
\right.\] 
There exists a sequence of polynomials $Q_k$, $k=1,2,...$ converging to $w$ uniformly on compact subsets of $ \ {{\rm \Omega}}_1\cup {{\rm \Omega}}_2$. Since ${{\rm \Omega}}_1\cup {{\rm \Omega}}_2$ is open, Weierstrass theorem applies. Thus,
\[{lim}_{k\to \infty }\left(D^{\left(\ell \right)}Q_k\right)=\left(D^{\left(\ell \right)}w\right)\] 
uniformly on every compact subset of ${{\rm \Omega}}_1\cup {{\rm \Omega}}_2$, whenever $\ell =\left({\ell }_1,\dots ,{\ell }_n\right)\in {{\mathbb N}}^n$. Therefore, there exists a $\ k_0\in {\mathbb N}$ such that, for all $\ell $ in the finite set $\mathfrak{F}$ we have

\begin{center}
 ${sup}_{z\in K}\left|Q\left(z\right)-P\left(z\right)\right|<\epsilon$\textit{ and }${sup}_{z\in \overline{{\rm \Omega}}\bigcap \overline{D\left(0;\delta\right)}}\left|D^{\left(\ell \right)}Q\left(z\right)\right|<\epsilon$
 \end{center}

\noindent \textit{ }where\textit{ }$Q\equiv Q_{k_0}$. \\

Combination of Lemma 3.2 with Theorem 3 in \cite{Bay}, gives the following result. \\
\\
\textbf{Theorem 3.3} \textit{Let }${\mathbb M}$\textit{ be an infinite subset of }${\mathbb N}$\textit{. Let }${\rm \Omega}\subset{{\mathbb C}}^n$\textit{ be an open convex set and }$\zeta\in \overline{{\rm \Omega}}$\textit{. Let also} ${\left(N_j\right)}_{j=0,1,2,\dots }$ \textit{be any enumeration of} ${{\mathbb N}}^n$\textit{. Then there exists a holomorphic function }$f\in A^{\infty }\left({\rm \Omega}\right)$\textit{ such that the partial sums} \newline $S_\lambda\left(f,\zeta\right)\left(z\right)=\sum^\lambda_{j=0}{a_{N_j}\left(f,\zeta\right){\left(z-\zeta\right)}^{N_j}}=\sum^\lambda_{j=0}{\left(\frac{1}{\left(N_{j,1}\right)!\dots \left(N_{j,n}\right)!}\frac{{\partial }^{N_{j,1}+\dots +N_{j,n}}f}{\partial z^{N_{j,1}}_1\dots \partial z^{N_{j,n}}_n}\right)\left(\zeta\right){\left(z-\zeta\right)}^{N_j}}$\textit{, }$\lambda=0,1,2,\dots $\textit{ satisfy the following.}

\textit{For every compact convex set} $K\subset \subset {{\mathbb C}}^n$\textit{ which is disjoint from }$\overline{{\rm \Omega}}$\textit{ and every holomorphic polynomial }$P\left(z\right)$\textit{, }$z\in{{\mathbb C}}^n$\textit{, there exists a sequence }$\left(\lambda_N\in{\mathbb M}\ :N=1,2,\dots \right)$\textit{ such that} 
\begin{center}
$\left[S_{\lambda_N}\left(f,\zeta\right)\right]\left(z\right)\rightarrow P\left(z\right)$\textit{ uniformly on }$K$\textit{, as }$N\rightarrow \infty$
\end{center}

\textit{and }
\begin{center}
$\left(D^{\left(\ell \right)}\left[S_{\lambda_N}\left(f,\zeta\right)\left(z\right)\right]\right)\rightarrow \left(D^{\left(\ell \right)}f\left(z\right)\right)$\textit{ uniformly on each compact subset of }$\overline{{\rm \Omega}}$\textit{, as }$N\rightarrow \infty$\textit{, for every}$\ \ell =\left({\ell }_1,\dots ,{\ell }_n\right)\in {{\mathbb N}}^n.$
\end{center}
\textit{Here }$S_\lambda\left(f,\zeta\right)\left(z\right)=\sum^\lambda_{j=0}{a_{N_j}\left(f,\zeta\right){\left(z-\zeta\right)}^{N_j}}$\textit{ with}
$$a_{N_j}\left(f,\zeta\right)=\frac{1}{\left(N_{j,1}\right)!\dots \left(N_{j,n}\right)!}\frac{{\partial }^{N_{j,1}+\dots +N_{j,n}}f}{\partial z^{N_{j,1}}_1\dots \partial z^{N_{j,n}}_n}\left(\zeta\right)$$
\textit{ and the }$a_{N_j}\left(f,\zeta\right)$\textit{ are well defined even for }$\zeta\in \partial {\rm \Omega}$\textit{, by the continuity of every }$D^{\left(N_j\right)}f$\textit{ on }$\overline{{\rm \Omega}}$\textit{. The set of all such functions }$f$\textit{ is dense and }${{\mathcal G}}_\delta$\textit{ in}$\ \ A^{\infty }\left({\rm \Omega}\right)$ \textit{and contains a dense vector subspace except }$0$.\\

Combining Theorem 3.3 with Theorem 3 in \cite{Bay}, we obtain the following. \\
\\
\textbf{Corollary 3.4} \textit{Let }${\mathbb M}$\textit{ be an infinite subset of }${\mathbb N}$\textit{. Let }${\rm \Omega}\subset {{\mathbb C}}^n$\textit{ be an open convex set. Let also }${\left(N_j\right)}_{j=0,1,2,\dots }$\textit{ be any enumeration of} ${{\mathbb N}}^n$\textit{. Then there exists a holomorphic function }$f\in A^{\infty }\left({\rm \Omega}\right)$\textit{ such that the following hold.}

\textit{For every compact convex set} $K\subset \subset {{\mathbb C}}^n$\textit{ disjoint from }$\overline{{\rm \Omega}}$\textit{ and every holomorphic polynomial }$P\left(z\right)$\textit{, }$\left(z\in{{\mathbb C}}^n\right)$\textit{, there exists a sequence }${\left(\lambda_N\in{\mathbb M}\right)}_{N=1,2,\dots }$\textit{ with the following properties.} \textit{For each compact set }$L\subset \overline{{\rm \Omega}}$\textit{, we have}
\begin{center}
${sup}_{\zeta\in L,\ z\in K}\left|S_{\lambda_N}\left(f,\zeta\right)\left(z\right)-P\left(z\right)\right|\to 0$\textit{ as }$N\to +\infty $
\end{center}
\textit{and}
\begin{center}
${sup}_{\zeta\in L,\ z\in L}\left|D^{\left(\ell \right)}\left[S_{\lambda_N}\left(f,\zeta\right)\right]\left(z\right)-D^{\left(\ell \right)}f\left(z\right)\right|\to 0$\textit{ as }$N\to +\infty $
\end{center}
\textit{whenever }$\ell =\left({\ell }_1,\dots ,{\ell }_n\right)\in {{\mathbb N}}^n$. \textit{Here} 
$$S_\lambda\left(f,\zeta\right)\left(z\right)=\sum^\lambda_{j=0}{a_{N_j}\left(f,\zeta\right){\left(z-\zeta\right)}^{N_j}}$$
\textit{with}
$$a_{N_j}\left(f,\zeta\right)=\frac{1}{\left(N_{j,1}\right)!\dots \ \left(N_{j,N}\right)!}\frac{{\partial }^{N_{j,1}+\dots +N_{j,n}}f}{\partial z^{N_{j,1}}_1\dots \partial z^{N_{j,n}}_n}\left(\zeta\right)$$
\textit{and the }$a_{N_j}\left(f,\zeta\right)$\textit{ are well defined even for }$\zeta\in \partial {\rm \Omega}$\textit{, by the continuity of every }$D^{\left(N_j\right)}f$\textit{ on }$\overline{{\rm \Omega}}$\textit{. The set of all such functions }$f$\textit{ is dense and }${{\mathcal G}}_\delta$\textit{ in}$\ \ A^{\infty }\left({\rm \Omega}\right)$ \textit{and contains a dense vector subspace except }$0$.

\section{Further Results}

 Since the disjoint union of a polynomially convex compact set with a finite number of smooth compact curves in its complement is polynomially convex (\cite{Sto}), our method yields the following\\.
\\
\noindent \textbf{Proposition 4.1} \textit{Let }${\rm \Omega}\subset {{\mathbb C}}^n$\textit{ be an open set such that }${\rm \Omega}={\overline{{\rm \Omega}}}^{{\rm o}}$\textit{ and }$\overline{{\rm \Omega}}$\textit{ is a polynomially convex compact set. Let }${{\mathbb P}}^{\infty }\left({\rm \Omega}\right)$\textit{ be the closure in }$A^{\infty }\left({\rm \Omega}\right)$\textit{ of the set of polynomials. Let also }${\mathbb M}\subset {\mathbb N}$\textit{ be an infinite set and }${\left(N_j\right)}_{j=0,1,2,\dots }$\textit{ be an enumeration of }${{\mathbb N}}^n$\textit{. Let also }${\left(K_m\right)}_{m=1,2,\dots }$\textit{ be a sequence of compact subsets of }${\overline{{\rm \Omega}}}^{{\rm c}}$\textit{, such that each }$K_m$\textit{ is a disjoint union of a finite number of smooth compact curves. Then, there exists a holomorphic function }$f\in {{\mathbb P}}^{\infty }\left({\rm \Omega}\right)$\textit{ such that, for every }$m\in \left\{1,2,\dots \right\}$\textit{ and every analytic polynomial P, there is a sequence }${\left(\lambda_N\in {\mathbb M}\right)}_{N=1,2,\dots }$\textit{ such that }

${sup}_{\zeta\in \overline{{\rm \Omega}},\ z\in K_m}\left|D^{\left(\ell \right)}\left[S_{\lambda_N}\left(f,\zeta\right)\right]-D^{\left(\ell \right)}P\left(z\right)\right|\to 0$\textit{ as }$N\to +\infty $\textit{\newline and}

${sup}_{\zeta\in \overline{{\rm \Omega}},\ z\in \overline{{\rm \Omega}}}\left|D^{\left(\ell \right)}\left[S_{\lambda_N}\left(f,\zeta\right)\right]\left(z\right)-D^{\left(\ell \right)}f\left(z\right)\right|\to 0$\textit{ as }$N\to +\infty $\textit{\newline for every derivation symbol }$\ell =\left({\ell }_1,\dots ,{\ell }_n\right)\in {{\mathbb N}}^n.$

\textit{The set of all such functions }$f\in {{\mathbb P}}^{\infty }\left({\rm \Omega}\right)$\textit{ is a dense and }${{\mathcal G}}_\delta$\textit{ subset of }${{\mathbb P}}^{\infty }\left({\rm \Omega}\right)$\textit{ and contains a vector space except }$0$\textit{.}\\

 To give another application of our method, let us recall a well known result by G. Henkin. \\
\\
\noindent \textbf{Lemma 4.2} (\cite{Hen}; see also Theorem 2.1 in page 280 of \cite{Ran}) \textit{Let }${\rm \Omega}\subset {{\mathbb C}}^n$\textit{ be a strictly pseudoconvex open domain with }$C^2$\textit{ boundary. Then the functions which are holomorphic in a neighborhood of }$\overline{{\rm \Omega}}$\textit{ are dense in the algebra }$A\left({\rm \Omega}\right)$\textit{ with respect to the supremum norm of }$\overline{{\rm \Omega}}$\textit{.}\\

In view of Lemma 4.2, our method yields immediately the following.\\
\\
\textbf{Corollary 4.3} \textit{Let }${\rm \Omega}\subset {{\mathbb C}}^n$\textit{ be a strictly pseudoconvex domain with }$C^2$\textit{ boundary. Suppose }$\overline{{\rm \Omega}}$\textit{ is a polynomially convex compact subset of }${{\mathbb C}}^n$\textit{. Let }${\mathbb M}\subset {\mathbb N}$\textit{ be an infinite set and }${\left(N_j\right)}_{j=0,1,2,\dots }$\textit{ be an enumeration of }${{\mathbb N}}^n$\textit{. Let also }${\left(K_m\right)}_{m=1,2,\dots }$\textit{ be a sequence of compact subsets of }${\overline{{\rm \Omega}}}^{{\rm c}}$\textit{, such that each }$K_m$\textit{ is a disjoint union of a finite number of smooth compact curves. Then, there exists a function }$f\in A\left({\rm \Omega}\right)$\textit{ such that, for every }$m\in \left\{1,2,\dots \right\}$\textit{ and every analytic polynomial P, there is a sequence }${\left(\lambda_N\in {\mathbb M}\right)}_{N=1,2,\dots }$\textit{ such that }\\

\textit{for every compact set }$L\subset {\rm \Omega}$\textit{ the following hold}

\begin{center}
 ${sup}_{\zeta\in L,\ z\in K_m}\left|D^{\left(\ell \right)}\left[S_{\lambda_N}\left(f,\zeta\right)\right]\left(z\right)-D^{\left(\ell \right)}P\left(z\right)\right|\to 0$\textit{ as }$N\to +\infty $\\
\textit{for  every  derivation  symbol}  $\ell =\left({\ell }_1,\dots ,{\ell }_n\right)\in {{\mathbb N}}^n$
\end{center}

\textit{and}

\begin{center}
${sup}_{\zeta\in L,\ z\in \overline{{\rm \Omega}}}\left|\left[S_{\lambda_N}\left(f,\zeta\right)\right]\left(z\right)-f\left(z\right)\right|\to 0$\textit{ as }$N\to +\infty $\textit{.} 
\end{center}

\textit{The set of all such functions }$f\in A\left({\rm \Omega}\right)$\textit{ is a dense and }${{\mathcal G}}_\delta$\textit{ subset of }$A\left({\rm \Omega}\right)$\textit{ and contains a vector space except }$0$\textit{.}\\
\\
\noindent \textbf{Proposition 4.4} \textit{Let }${\rm \Omega}\subset {{\mathbb C}}^n$\textit{ be an open set which admits an exhausting family of compact sets, which are polynomially convex. Let }${\mathbb M}\subset {\mathbb N}$\textit{ be an infinite set and }${\left(N_j\right)}_{j=0,1,2,\dots }$\textit{ be an enumeration of }${{\mathbb N}}^n$\textit{. Let also }${\left(K_m\right)}_{m=1,2,\dots }$\textit{ be a sequence of compact subsets of }${{\rm \Omega}}^{{\rm c}}$\textit{, such that each }$K_m$\textit{ is a disjoint union of a finite number of smooth compact curves. Then, there exists a holomorphic function }$f\in {\mathcal O}\left({\rm \Omega}\right)$\textit{ such that, for every }$m\in \left\{1,2,\dots \right\}$\textit{ and every analytic polynomial P, there is a sequence }${\left(\lambda_N\in {\mathbb M}\right)}_{N=1,2,\dots }$\textit{ such that}\\

\textit{for every compact set }$L\subset {\rm \Omega}$\textit{ the following hold}

\begin{center}
${sup}_{\zeta\in L,\ z\in K_m}\left|D^{\left(\ell \right)}\left[S_{\lambda_N}\left(f,\zeta\right)\right]\left(z\right)-D^{\left(\ell \right)}P\left(z\right)\right|\to 0$\textit{ as }$N\to +\infty $
\end{center}

\textit{and }

\begin{center}
${sup}_{\zeta\in L,\ z\in L}\left|D^{\left(\ell \right)}\left[S_{\lambda_N}\left(f,\zeta\right)\right]\left(z\right)-D^{\left(\ell \right)}f\left(z\right)\right|\to 0$\textit{ as }$N\to +\infty $\textit{\newline for every derivation symbol }$\ell =\left({\ell }_1,\dots ,{\ell }_n\right)\in {{\mathbb N}}^n$\textit{.} 
\end{center}

\textit{The set of all such functions }$f\in {\mathcal O}\left({\rm \Omega}\right)$\textit{ is a dense and }${{\mathcal G}}_\delta$\textit{ subset of }${\mathcal O}\left({\rm \Omega}\right)$\textit{ endowed with the topology of uniform convergence on compact subsets of }${\rm \Omega}$\textit{ and contains a vector space except }$0$\textit{.}\\

\noindent \textbf{Corollary 4.5} \textit{Let }${\rm \Omega}\subset {{\mathbb C}}^n$\textit{ be an open set. Suppose }$\varrho :{\rm \Omega}\to {\mathbb C}$\textit{ is a continuous plurisubharmonic exhaustion function for }${\rm \Omega}$\textit{ so that each sublevel set }${{\rm \Omega}}_c=\left\{z\in {\rm \Omega}:\varrho \left(z\right)<c\right\}$\textit{ is relatively compact in }${\rm \Omega}$\textit{. Let }${\mathbb M}\subset {\mathbb N}$\textit{ be an infinite set and }${\left(N_j\right)}_{j=0,1,2,\dots }$\textit{ be an enumeration of }${{\mathbb N}}^n$\textit{. Let also }${\left(K_m\right)}_{m=1,2,\dots }$\textit{ be a sequence of compact subsets of }${{\rm \Omega}}^{{\rm c}}$\textit{, such that each }$K_m$\textit{ is a disjoint union of a finite number of smooth compact curves. Then, there exists a holomorphic function }$f\in {\mathcal O}\left({\rm \Omega}\right)$\textit{ such that, for every }$m\in \left\{1,2,\dots \right\}$\textit{ and every analytic polynomial P, there is a sequence }${\left(\lambda_N\in {\mathbb M}\right)}_{N=1,2,\dots }$\textit{ such that}\\

\textit{for every compact set }$L\subset {\rm \Omega}$\textit{ the following hold}\\
\begin{center}
${sup}_{\zeta\in L,\ z\in K_m}\left|D^{\left(\ell \right)}\left[S_{\lambda_N}\left(f,\zeta\right)\right]\left(z\right)-D^{\left(\ell \right)}P\left(z\right)\right|\to 0$\textit{ as }$N\to +\infty $
\end{center}

\textit{and}

\begin{center}
${sup}_{\zeta\in L,\ z\in L}\left|D^{\left(\ell \right)}\left[S_{\lambda_N}\left(f,\zeta\right)\right]\left(z\right)-D^{\left(\ell \right)}f\left(z\right)\right|\to 0$\textit{ as }$N\to +\infty $\textit{\newline for every derivation symbol }$\ell =\left({\ell }_1,\dots ,{\ell }_n\right)\in {{\mathbb N}}^n$\textit{.}
\end{center}

\textit{The set of all such functions }$f\in {\mathcal O}\left({\rm \Omega}\right)$\textit{ is a dense and }${{\mathcal G}}_\delta$\textit{ subset of }${\mathcal O}\left({\rm \Omega}\right)$\textit{ endowed with the topology of uniform convergence on compact subsets of }${\rm \Omega}$\textit{ and contains a vector space except }$0$\textit{. }

\noindent \textit{Proof} It is well known that each region ${{\rm \Omega}}_c$ is a Runge domain (see, for instance, Theorem 1.3.7, page 25, in \cite{Stou}). Since every Runge domain can be exhausted by polynomially convex compact sets, the open set ${\rm \Omega}\subset {{\mathbb C}}^n$ admits an exhausting family of compact sets, which are polynomially convex. Application of Proposition 4.4 completes the proof.


\begin{thebibliography}{}
%
%
\bibitem{Back} U. Backlund and A. $F\ddot{a}llstr\ddot{o}m$, \textit{The polynomial hull of unions of convex sets in} ${{\mathbb C}}^n$, Colloquim Mathematicum, LXX (1996)(1) 7-11

\bibitem{Bay}  F. Bayart, K.G. Grosse-Erdmann, V. Nestoridis and C. Papadimitropoulos: \textit{Abstract Theory of Universal Series and Applications}, Proc. London Math. Soc., Volume 96, Issue 3, 2008, pages 417-463 

\bibitem{Bed}  E. Bedford and S. Pinchuk, \textit{Domains in }${{\mathbb C}}^{n+1}$ \textit{with noncompact automorphism group}, J. Geom. Anal. 1 (1991), 165--192

\bibitem{Chui} C. Chui and M.N. Parnes: \textit{Approximation by overconvergence of power series}, Journal of Mathematical Analysis and Applications, Volume 36, 1971, pages 693-696 

\bibitem{Clo}  R. $Clou\hat{a}tre$: \textit{Universal power series in} ${{\mathbb C}}^N$, Canad. J. Math. Volume 54, Number 2, 2011 pages 230-236

\bibitem{For}  J. E. Fornaess  and N. Sibony, \textit{Complex dynamics in higher dimension}\textbf{, }Proceedings Berkeley, in \textbf{Several Complex Variables}, M.Schneider,Y.T. Siu editors, MSRI publications, (1999), 273-296

\bibitem{Forst}  F. Forstneri$\check{c}$: \textbf{\textit{Stein Manifolds and Holomorphic Mappings. The Homotopy Principle in Complex Analysis}}, A Series of Modern Surveys in Mathematics 56, Ergebnisse der Mathematik und ihrer Grenzgebiete, 3. Folge, Springer, 2011

\bibitem{Gros1}  K.-G. Grosse-Erdmann: \textit{Holomorphe Monster und universelle Funktionen}, Mitt. Math. Sem. Giessen, Volume 176, 1987, iv+84 pages 

\bibitem{Gros2}  K.-G. Grosse-Erdmann: \textit{Universal families and hypercyclic operators}, Bull. Amer. Math. Soc., Volume 36, 1999, pages 345-381 

\bibitem{Hen}  G. Henkin: \textit{Integral representations of functions holomorphic in strictly pseudoconvex domains and some applications}, Mat. Sbornik, Vol. 78, 1969, pages 611-632 (in Russian)

\bibitem{Kah} J.-P. Kahane: \textit{Baire's Category theorem and trigonometric series}, J. Anal. Math., Volume 80, 2000, pages 143-182

\bibitem{Kal}  E. Kallin: \textit{Polynomial convexity}: \textit{The three spheres problem}, in: Proceedings of the Conference on Complex Analysis, Minneapolis 1964, H. Röhrl, A. Aeppli and E. Calabi (eds.), Springer, 1965, pages 301--304

\bibitem{Kal2}  E.Kallin: \textit{Fat polynomially convex sets}, in \textbf{Function Algebras} (Chicago: Scott Foresman)(1966), pp. 149--152\textbf{; Proceedings of International Symposium on Function Algebras}, Tulane University, 1965

\bibitem{Khu}  G. Khud$\check{a}$iberganov, \textit{On polynomial and rational convexity of unions of compacts in }${{\mathbb C}}^n$, Izv. Vuz. Mat. 2 (1987), 70--74 (in Russian)

\bibitem{Kyt}  A. M. Kytmanov and G. Khud$\check{a}$iberganov, \textit{An example of a non-polynomially convex compact set consisting of three non-intersecting ellipsoids}, Sibirsk. Mat. Zh. 25(5) (1984), 196--198 (in Russian)

\bibitem{Luh}  W. Luh: \textit{Approximation analytischer Funktionen durch uberkonvergente Potenzreihen und deren Matrix-Transformierten}, Mitt. Math. Sem. Giessen, Volume 88, 1970, pages 1-56 

\bibitem{Mel1} A. Melas and V. Nestoridis: \textit{Universality of Taylor series as a generic property of holomorphic functions}, Adv. Math., Volume 157, 2001, Pages 138-176

\bibitem{Mel2}  A. Melas and V. Nestoridis: \textit{On various types of universal Taylor series}, Complex Variables Theory, Volume 44, Issue 3, 2001, pages 245-258

\bibitem{Mul}  J. M$\ddot{u}$ller, V. Vlachou and A. Yavrian: \textit{Universal overconvergence and ostrowski-gaps}, Bull. Lomdon Math. Soc., Volume 38, Issue 4, 2006, pages 597-606

\bibitem{Nest1}  V. Nestoridis: \textit{Universal Taylor series}, Annales de l'Institut Fourier, Volume 46, 1996, pages 1293-1306

\bibitem{Nest2}  V. Nestoridis: \textit{An extension of the notion of universal Taylor series}, in N. Papamichael, S. Ruscheweyh, E. B, Saff (eds) Proceedings of the 3rd CMFT Conference on Computational Methods and Function Theory, 1997, Nicosia, Cyprus, October 13-17, 1997, World Scientific Ser. Approx. Decompos. 11(1999), page 421-430

\bibitem{Nest3}  V. Nestoridis: \textit{A strong notion of universal Taylor series}, J. London Math. Soc., Volume 68, Issue 2, 2003, pages 712-724

\bibitem{Nest4}  V. Nestoridis and Chr. Papadimitropoulos: \textit{Abstract theory of universal series and an application to Dirichlet Series}, C. R. Acad. Sci. Paris, Ser. I 341, 2005, pages 539-543

\bibitem{Pae}   de Paepe, P.J.: \textit{Eva Kallin's lemma on polynomial convexity}, The Bulletin of the London Mathematical Society, 33(2001), 1--10

\bibitem{Pet} H. Peters and E. Fornaess Wold, \textit{Non-autonomous basins of attraction and their boundaries}, Journal of Geometric Analysis, \textit{15}\eqref{GrindEQ__1_} (2005), 123- 135

\bibitem{Ran}  M. Range: \textbf{\textit{Holomorphic Functions and Integral Representations in Several Complex Variables}}, Springer-Verlag, New York 1986. 2nd corrected printing 1998

\bibitem{Ros}  J.-P. Rosay, \textit{The polynomial hull of non-connected tube domains, and an example of E. Kallin}, Bull. London Math. Soc. 21 (1989), 73--78

\bibitem{Sto}  G. Stolzenberg: \textit{Uniform approximation on smooth curves}, Acta Math. 115 (1966), 185-198

\bibitem{Stou}  E. L. Stout, \textbf{\textit{Polynomial Convexity}}, Series: \textbf{Progress in Mathematics}, Vol. 261, Boston, Mass. : Birkh$\ddot{a}$user; London: Springer [distributor], 2007

\bibitem{Wer}  J. Wermer, \textit{Polynomial approximation on an arc in }${{\mathbb C}}^3$, Ann. of Math. 62 (1955), 269--270

\bibitem{Ex}   ---, \textit{An example concerning polynomial convexity}, Math. Ann. 139 (1959), 147--150
\end{thebibliography}
\end{document}